\documentclass[11pt,a4]{article}
\input{amssym}
\title{{\Large \bf The Annihilating-Ideal  Graph  of Commutative  Rings I\footnote{ The research
 of the first author was in part supported by
a grant from IPM (No. 87160026).}}}
\author{{\large\bf M. Behboodi$^{a,b}$\footnote{Corresponding
author.} and Z. Rakeei$^a$ }
\\{\footnotesize $^{a}$Department of Mathematical Science,
 Isfahan University of Technology,  } \\{\footnotesize P. O. Box: 84156-83111, Isfahan, Iran }\\{\footnotesize $^{b}$School of Mathematics,
  Institute for Research in Fundamental Sciences (IPM), } \\{\footnotesize P. O. Box: 19395-5746, Tehran, Iran}\\{\footnotesize
mbehbood@cc.iut.ac.ir, \ \ sanamzhr@yahoo.com}}

\def\be{\begin{enumerate}}
\def\ee{\end{enumerate}}

\newtheorem{ttheo}{Theorem}[section]
\newtheorem{ccoro}[ttheo]{Corollary}
\newtheorem{llem}[ttheo]{Lemma}
\newtheorem{eexam}[ttheo]{Example}
\newtheorem{rrem}[ttheo]{Remark}
\newtheorem{ppro}[ttheo]{Proposition}

\newenvironment{pproof}{\noindent{\bf Proof. }}{}

\date{}
 \begin{document}
  \maketitle
\begin{center}
  {\small\bf Abstract}
  \end{center}

  {\small Let $R$ be  a commutative ring with ${\Bbb{A}}(R)$ its set of ideals  with nonzero annihilator.
   In this paper and its sequel, we introduce and investigate the {\it annihilating-ideal graph}  of $R$, denoted by  ${\Bbb{AG}}(R)$.
   It is the (undirected) graph with  vertices  ${\Bbb{A}}(R)^*:={\Bbb{A}}(R)\setminus\{(0)\}$, and
   two  distinct  vertices  $I$ and $J$ are adjacent if and only if  $IJ=(0)$.  First, we study some
     finiteness conditions of ${\Bbb{AG}}(R)$. For instance,  it is shown that if  $R$  is not a domain, then  ${\Bbb{AG}}(R)$ has
    ACC   (resp., DCC) on vertices  if and only if
 $R$ is  Noetherian  (resp.,  Artinian). Moreover,  the set of vertices of  ${\Bbb{AG}}(R)$ and the set of nonzero proper ideals of
  $R$ have the same cardinality when $R$ is either an  Artinian or a  decomposable ring. This yields for a ring $R$,
  ${\Bbb{AG}}(R)$ has $n$ vertices  $(n\geq 1)$ if  and only if $R$ has only $n$  nonzero proper
 ideals. Next, we study  the connectivity  of ${\Bbb{AG}}(R)$. It is shown that ${\Bbb{AG}}(R)$ is a  connected  graph and
$diam(\Bbb{AG})(R)\leq 3$ and if ${\Bbb{AG}}(R)$ contains a cycle,
then $gr({\Bbb{AG}}(R))\leq 4$. Also, rings $R$ for which the
graph ${\Bbb{AG}}(R)$ is complete or star, are characterized, as
well as rings $R$ for which every vertex of ${\Bbb{AG}}(R)$ is a
prime (or maximal) ideal.  In Part II we shall study the diameter
and
coloring of annihilating-ideal graphs.}\vspace{2mm}\vspace{2mm} \vspace{4mm}\\

 \noindent{\bf 0. Introduction}
 \vspace{2mm}

 In the literature, there are many papers on assigning a
graph to a  ring, a group, semigroup or a  module (see for example
[1-11,13-17]). In fact, the concept of the zero divisor graph of a
commutative ring $R$ was first introduced by Beck [6], where he
was mainly interested in colorings. In his work all elements of
the ring were vertices of the graph. This investigation of
colorings of a commutative ring was then continued by Anderson and
Naseer in [2]. Let $Z(R)$ be the set of zero-divisors of $R$. In
[3], Anderson and Livingston associate a graph, $\Gamma(R)$, to
$R$ with vertices $Z(R)\setminus \{0\}$, the set of nonzero
zero-divisors of $R$, and for distinct $x$, $y\in
Z(R)\setminus\{0\}$, the vertices $x$ and $y$ are adjacent if and
only if $xy=0$. In [18], Sharma and Bhatwadekar define another
graph on $R$, $G(R)$, with vertices as elements of $R$, where two
distinct vertices $a$ and $b$ are adjacent if and only if $Ra + Rb
= R$ (see also [15], in which, the notion of  ``comaximal graph of
commutative rings" is investigated). Recently Anderson and Badawi
in [4] have  introduced  and investigated  the total graph of $R$,
denoted by $T(\Gamma(R))$. It is the (undirected) graph with all
elements of $R$ as vertices, and for distinct $x$, $y\in R$, the
vertices $x$ and $y$ are adjacent if and only if $x+y\in Z(R)$.

In ring theory, the structure of a ring $R$ is  closely tied to
ideal's behavior more than elements', and so it is deserving to
define a graph whose vertices are ideals instead of elements.

Throughout this paper  $R$ will be a commutative ring with
identity. We denote the set of all proper ideals of $R$ by
$\Bbb{I}(R)$. We name an ideal $I$ of $R$, an {\it
annihilating-ideal} if there exists a nonzero ideal $J$ of $R$
such that  $IJ=(0)$, and use the notation $\Bbb{A}(R)$ for the set
of all annihilating-ideals of $R$. Also  $Spec(R)$, $Max(R)$ and
$J(R)$ denote the sets of prime ideals,  maximal ideals and the
Jacobson radical of $R$, respectively.

    In this paper,  we  define {\it annihilating-ideal graph} of  $R$, denoted by  ${\Bbb{AG}}(R)$,  as  a  graph   with vertices
    $\Bbb{A}(R)^* =\Bbb{A}(R)\hspace{-1mm}\setminus\{(0)\}$,
   where distinct  vertices  $I$ and $J$ are adjacent if and only if  $IJ=(0)$. We investigate  the interplay
       between the  graph-theoretic
   properties of ${\Bbb{AG}}(R)$ and  the ring-theoretic properties of $R$.  In Section 1, we study some finiteness conditions of annihilating-ideal graphs.
   For instance,  it is shown that if  $R$  is not a domain, then  ${\Bbb{AG}}(R)$ has
    assenting chain condition ($ACC$) (resp., descending chain condition   ($DCC$)) on vertices  if and only if
 $R$ is a Noetherian  (resp., an Artinian) ring.  Also, it is shown that  ${\Bbb{AG}}(R)$ has $n$ vertices
$(n\geq 1)$ if  and only if $R$ has only $n$  nonzero proper
ideals. These facts motivates us to the following  conjecture: for
a non-domain   ring $R$,  the set of vertices of ${\Bbb{AG}}(R)$
and the set of  nonzero proper
 ideals of $R$ have the same cardinality. The conjecture above is true for all Artinian rings as well as
 all decomposable rings (see Proposition 1.3 and Proposition 1.6).\\
   Recall that a graph $G$ is connected if there is a
path between every  two distinct vertices. For distinct vertices
$x$ and $y$ of $G$, let $d(x, y)$ be the length of the shortest
path from $x$ to $y$ and if there is no such path we define
$d(x,y)=\infty$. The  diameter of $G$ is $diam(G) = sup \{d(x,y) \
: \ x \ and \ y \ are \ distinct \ vertices \ of \ G\}$. The girth
of $G$, denoted by $gr(G)$, is defined as the length of the
shortest cycle in $G$ and $gr(G)=\infty$ if $G$ contains no
cycles. A graph in which each pair of distinct vertices is joined
by an edge is called a complete graph. Also, if a  graph $G$
contains one vertex to which all other vertices are joined and $G$
has no other edges, is called a star graph. In Section 2, the
connectivity  of the annihilating-ideal graphs are studied. It is
shown that for every  ring $R$, ${\Bbb{AG}}(R)$ is a connected
graph and $diam(\Bbb{AG})(R)\leq 3$,  and if ${\Bbb{AG}}(R)$
contains a cycle, then $gr({\Bbb{AG}}(R))\leq 4$. Also, rings $R$
for which the graph ${\Bbb{AG}}(R)$ is a complete (or  star) graph
are characterized, as well as rings $R$ for which every vertex of
${\Bbb{AG}}(R)$ is a prime (or maximal) ideal. In part II we shall continue the study of this construction via diameter and coloring.\\\\

%\newpage
\noindent{\bf 1. Finiteness conditions of annihilating-ideal graphs}\\

 For a  ring $R$,  $soc(R)$ is the sum of all minimal
ideals of $R$ (if there are no minimal  ideals, this sum is
defined to be zero). Also, if $X$ is an element or a subset of a
ring  $R$, we define the annihilator of $X$ in $R$  by
$Ann(X)=\{r\in R \ | \ rX=(0)\}$. An ideal $I$ of a ring  $R$ is
called an annihilator
ideal if $I=Ann(x)$ for some $x\in R$. Also, for a set $A$ we denote the cardinal number of $A$ by $|A|$. \\

   Let $R$ be a   ring. We say that the annihilating-ideal graph ${\Bbb{AG}}(R)$ has $ACC$ (resp., $DCC$) on vertices   if
$R$ has $ACC$ (resp., $DCC$) on ${\Bbb{A}}^*(R)$. \\

\noindent{\bf Theorem 1.1.} {\it Let $R$  be a
 non-domain ring. Then  ${\Bbb{AG}}(R)$ has  $ACC$ (resp., $DCC$) on vertices  if and only if
 $R$ is a Noetherian  (resp., an Artinian) ring.}\\

\noindent{\bf  Proof.}   Suppose that ${\Bbb{AG}}(R)$ has $ACC$
(resp., $DCC$) on vertices. Let $0\neq x\in Z(R)$ and $P=Ann(x)$.
Then
$$\{I~:~I\trianglelefteq R,~I\subseteq Rx\}\bigcup\{I~:~I\trianglelefteq R,~I\subseteq P\}\subseteq
\Bbb{A}(R).$$ It follows that  the $R$-modules $Rx$ and $P$ have
$ACC$ (resp., $DCC$) on submodules i.e., $Rx$ and $P$ are
Noetherian  (resp., Artinian) $R$-modules. Since $Rx\cong R/P$,
by [19, Proposition 7.17], $R$ is a Noetherian (resp., an Artinian) ring.  The converse is clear. $\square$\\

In view of above theorem and Cohen's theorem, one may naturally
ask, when $Spec(R)\cap\Bbb{A}(R)^*\neq \emptyset$ and every prime
ideal $P\in \Bbb{A}(R)^*$ is finitely generated, is $R$
Noetherian? The answer is no! The
following example gives a non-Noetherian ring $R$ for which every  prime vertex of $\Bbb{AG}(R)$ is finitely generated.\\

\noindent{\bf Example 1.2.} Let $K$ be a field and $D:=K[\{x_i :
i\in \Bbb{Z}\}]$. The domain $D$ is a unique factorization domain
and $Q=x_1D$ is a principal height one prime of $D$. Let
$R:=D/x_1^2D$  and $P_0:=Q/x_1^2D$. Then  $P_0$ is a prime ideal
of $R$ with $P_0^2=(0)$. Thus $P_0\in {\Bbb{AG}}(R)$ and
$P_0\subseteq P$ for every $P\in Spec(R)$. Moreover, $x_1^2D$ is a
$Q$-primary ideal of $D$. It follows that the zero ideal of $R$ is
$P_0$-primary. Hence $IJ=(0)$ in $R$ with $I\neq(0)$ and
$J\neq(0)$, implies both $I$ and $J$ are contained in
$P_0=Q/x_1^2D$. Therefore,  $P_0=Z(R)$ is a
  cyclic ideal  and
$Spec(R)\cap\Bbb{A}(R)^*=\{P_0\}$,  but $R$ is not a Noetherian ring.\\

The following proposition shows that  rings $R$ for which every
nonzero proper ideal $I$ of $R$ is a vertex of ${\Bbb{AG}}(R)$ are
  abundant.\\

 \noindent{\bf Proposition  1.3.} {\it Let $R$ be an Artinian
  ring. Then  every nonzero proper ideal $I$ of $R$ is a
vertex of ${\Bbb{AG}}(R)$. }\\

\noindent{\bf  Proof.} Let $R$ be an Artinian   ring. Then  each
maximal ideal has a nonzero annihilator (since each is also a
minimal prime, and a minimal prime in a Noetherian ring always has
a nonzero annihilator). Thus every nonzero proper ideal $I$ of $R$
is a vertex of ${\Bbb{AG}}(R)$. $\square$\\

The following result shows that the set of all rings for which the
annihilating-ideal graphs are finite coincides with the set of all
(Artinian) rings with only finitely many  proper ideals. Moreover,
for a
  ring $R$, the graph ${\Bbb{AG}}(R)$ has $n$ vertices
$(n\geq 1)$ if and only if $R$ has only $n$  nonzero proper ideals.\\

\noindent{\bf Theorem   1.4.} {\it Let $R$ be a   ring. Then the following statements are equivalent.}\vspace{3mm}\\
(1) {\it  ${\Bbb{AG}}(R)$  is a finite graph.}\\
(2) {\it  $R$ has only finitely many  ideals.}\\
(3)  {\it Every vertex  of ${\Bbb{AG}}(R)$ has finite degree.}\vspace{2mm}\\
 {\it Moreover, ${\Bbb{AG}}(R)$ has $n$ $(n\geq 1)$
vertices if and
 only if $R$ has only $n$ nonzero  proper ideals.}\\

  \noindent{\bf Proof.} $(1)\Rightarrow(2)$.  Suppose that ${\Bbb{AG}}(R)$ is a finite
graph with $n$ ($n\geq 1$) vertices. Then by Theorem 1.1, $R$ is
an Artinian ring. Now by Proposition 1.3, every nonzero proper
ideal $I$ of $R$ is a vertex
of ${\Bbb{AG}}(R)$. Thus $|\Bbb{I}(R)^*|=n$ i.e., $R$ has only $n$  nonzero proper ideals.\\
$(2)\Rightarrow(3)$. Clear.\\
 $(3)\Rightarrow(1)$. Suppose that every vertex of ${\Bbb{AG}}(R)$ has finite degree, but ${\Bbb{AG}}(R)$ is an infinite
graph.  Let $I=Rx$ be a vertex of ${\Bbb{AG}}(R)$ and $J=Ann(I)$.
  If the set of $R$-submodules of $I$ (resp., $J$) is infinite, then $J$ (resp., $I$) has infinite degree, a contradiction.
  Thus the set of $R$-submodules of $I$ and  $J$ are finite and hence  $I$ and $J$ are Artinian $R$-module. Now, since
  $I\cong R/J$, $R$ is an Artinian ring. Thus by Proposition 1.3, every nonzero proper ideal
  of $R$ is a vertex of ${\Bbb{AG}}(R)$.  We  consider the following 2
  cases:\\
  \noindent{Case 1}: $R$ is a local ring. Since $I$ has finite
  number of $R$-submodules,  $R$ has a minimal ideal,  say $Rx_0$. Let $P=Ann(x_0)$. Then $P$ is the maximal ideal of $R$ and so
   every proper ideal of  $R$ is contained in $P$. This implies that  $Rx_0$ is adjacent to all other
    vertices  of ${\Bbb{AG}}(R)$. Since $Rx_0$ has finite degree, ${\Bbb{AG}}(R)$ is a finite graph. \\
   \noindent{Case 2}: $R$ is not local. Suppose that  $R=R_1\times R_2$ where $R_1$ and $R_2$ are   nonzero
  rings. Then for each ideal $I$ of $R_1$ the vertex $I\times (0)$ of ${\Bbb{AG}}(R)$ is adjacent to $(0)\times R_2$. Thus the set of ideals of $R_1$ is
  finite. Similarly, the set of ideals of $R_2$ is finite. This implies that  the set of ideals of $R$ is finite, i.e.,
   ${\Bbb{AG}}(R)$ is a finite graph.\\
   Finally, by Proposition 1.3, we conclude that ${\Bbb{AG}}(R)$ has $n$ $(n\geq 1)$
vertices if and  only if $R$ has only $n$ nonzero  proper ideals. $\square$\\

We have not found any examples of a non-domain    ring $R$ such
that  $|{\Bbb{A}}^*(R)|< |{\Bbb{I}}^*(R)|$. The lack of such
counterexamples, together with the fact that
$|{\Bbb{A}}^*(R)|=|{\Bbb{I}}^*(R)|$ where $R$ is an  Artinian ring
(see Proposition 1.3 and Theorem 1.4) or $R$ is a decomposable
ring (see Proposition
1.6),  motivates the following fundamental conjecture:\\

\noindent {\bf Conjecture 1.5.} {\it Let $R$ be a  non-domain
  ring. Then the set of vertices of ${\Bbb{AG}}(R)$ and
the set of  nonzero proper
 ideals of $R$ have the same cardinality.} \\

\noindent {\bf Proposition 1.6.} {\it The Conjecture 1.5 is true
for each decomposable    ring  $R$.} \\

\noindent {\bf Proof.} Suppose that $R=R_1\times R_2$. If both
$R_1$ and $R_2$ are Artinian i.e., $R$ is an Artinian ring, then
by Proposition 1.3, the proof is complete. Thus we can assume that
$R_1$ is not Artinian and so the set of nonzero proper
 ideals of $R_1$ is infinite. Clearly, every ideal of $R$ is of the form $I_1\times I_2$ where $I_1$ and
$I_2$ are ideals of $R_1$ and $R_2$, respectively. Note that $|
\Bbb{I}(R_1)\times \Bbb{I}(R_2)|$ equals $max \{ | \Bbb{I}(R_1)|,
| \Bbb{I}(R_2)| \}$. Also every ideal of the form $I\times (0)$ or
 $(0)\times J$, where $I$ is an ideal of $R_1$ and $J$ is an
ideal of $R_2$, is a vertex of  ${\Bbb{AG}}(R)$. Clearly
$|\{I\times (0) : I\in \Bbb{I}(R_1)\}|=|\Bbb{I}(R_1)|$  and
$|\{(0)\times J : J\in \Bbb{I}(R_2) \}|=|\Bbb{I}(R_2)|$ and hence
  $$|\Bbb{A}(R)|\geq max \{ | \Bbb{I}(R_1)|, | \Bbb{I}(R_2)| \}=|\Bbb{I}(R_1)|\times | \Bbb{I}(R_2)|=|\Bbb{I}(R)|.$$
  On the other hand $|\Bbb{A}(R)|\leq |{\Bbb{I}(R)}|$. It follows that  $|\Bbb{A}(R)^*|=|{\Bbb{I}(R)}^*|$. $\square$\\

 The following  result shows that if every nonzero proper ideal of a Noetherian ring
   $R$ is a vertex of
  ${\Bbb{AG}}(R)$, then $R$ is a  semilocal ring (i.e., $R$ has only finitely many maximal
  ideals).\\

 \noindent{\bf Proposition  1.7.} {\it Let $R$ be a Noetherian ring. If  all  nonzero proper ideals  of $R$ are
vertices  of ${\Bbb{AG}}(R)$, then $R$ is a semilocal ring. }\\

\noindent{\bf  Proof.}  Suppose that $\{P_n \ | \ n\in\Bbb{N}\}$
are distinct maximal ideals of $R$. Then by our assumption,
$Ann(P_n)\neq (0)$ for each $n\geq 1$. It follows that for each
$n$ there exists $x_n\in R$ such that $P_n=Ann(x_n)$. Thus for
each $n$, $x_n$ is in every maximal ideal $P_k$ for $k\neq n$. By
ACC, the chain $Rx_1\subsetneqq
 Rx_1+Rx_2\subsetneqq\cdots\subsetneqq Rx_1+\cdots+Rx_n\subsetneqq\cdots$ must stabilize, and each step
is proper since $P_1\cap P_2\cap\cdots\cap P_k$ is the annihilator
of $x_1R+x_2R+ \cdots + x_kR$ for each $k$. Therefore, $R$ has only finitely many maximal ideals, i.e., $R$ is semilocal.  $\square$\\

The following result shows that for each  Noetherian  ring $R$, if
 $R$ is not a domain, then at least one of the
vertices of $\Bbb{AG}(R)$ is a prime ideal.\\

  \noindent{\bf Proposition  1.8.} {\it Let  $R$ be a   Noetherian ring. Then either $\Bbb{A}(R)^*=\emptyset$ or  ${Spec(R)\cap\Bbb{A}(R)^*\neq\emptyset}$.}\\

\noindent{\bf  Proof.} Assume $\Bbb{A}(R)^*\neq\emptyset$ i.e.,
$R$ is not a domain. Clearly, the set of all annihilators of
nonzero elements of $R$  is a subset of $\Bbb{A}(R)^*$. On the
other hand by our hypothesis there exists an annihilator ideal $P$
of  $R$ which is maximal among all annihilators of nonzero
elements of $R$. Now  by [13, Theorem 6], $P$ is a prime ideal.
Since $R$ is not a domain,  $P\neq (0)$ and so
 $P\in Spec(R)\cap \Bbb{A}(R)^*$. $\square$ \\

  The following example shows that the  Noetherian hypothesis is needed in Proposition 1.8. \\

 \noindent{\bf Example 1.9.} Let $K$ be a field, $D:=K[\{x_i : i\in\Bbb{N}\}]$ (a unique factorization domain) and
$R=K[\{x_i : i\in \Bbb{N}\}]/(\{x_i^2 : i\in
 \Bbb{N}\})$. Let $\bar{x}_i=x_i+(\{x_i^2 : i\in
 \Bbb{N}\})$ for each $i\in\Bbb{N}$ and  $M=(\{\bar{x}_i :  i\in
 \Bbb{N}\})$. The ideal $M$ is simply the image of the maximal ideal $N=(\{x_i : i\in
 \Bbb{N}\})$ of $D$.
Moreover $N$ is the radical of $J=(\{x_i^2 : i\in
 \Bbb{N}\})$. Thus $R$
is a local zero-dimensional ring (but is not Noetherian). For each
$x_i$, the conductor of $x_i$ into $J$ is the height one
(principal) prime $x_iD$ (since $x_i$ is an irreducible=prime
element of $D$). Thus viewed as an element of $R$, the annihilator
of $x_i$ is $x_iR$. Obviously $\bigcap x_iD=(0)$ and this implies
$\bigcap x_iR=(0)$. Hence $Ann(M)=(0)$.
 Thus  $Spec(R)\cap\Bbb{A}(R)^*=\emptyset$.  $\square$ \\

We need the following two lemmas.\\

 \noindent{\bf Lemma 1.10.} {\it Let
$R$ be a   ring and $P$ be a maximal ideal such that $P\in
\Bbb{A}(R)^*$. Then
$P=Ann(x)$ for some  $0\neq x\in R$.}\\

\noindent{\bf  Proof.} Let $P\in \Bbb{A}(R)^*$ be a maximal ideal
of $R$. Then there exists $I\in \Bbb{A}^*(R)$,
  such that $IP=(0)$. Let $0\neq a\in I$. Then $aP=(0)$ and so $P\subseteq Ann(a)$.  Since $P$ is maximal,
  $P=Ann(a)$.  $\square$ \\

\noindent{\bf Lemma 1.11.} {\it Let $R$ be a   ring such that $R$
is not a field. Then for each minimal ideal $I$ of $R$,
$I$ and $Ann(I)$ are vertices  of ${\Bbb{AG}}(R)$.}\\

\noindent{\bf  Proof.} Suppose that $I$ is a minimal ideal of $R$.
Since $R$ is not a field, $I\neq R$. Now  by Brauer's Lemma (see
[13, 10.22]), either $I^2=(0)$ or $I=Re$ for some idempotent
$1\neq e\in R$. If $I^2=(0)$, then $I$ is a vertex of
${\Bbb{AG}}(R)$. If $I^2=Re$, then $R=Re\oplus R(1-e)$.
Clearly $ReR(1-e)=(0)$ and so $Re$ is a  vertex of ${\Bbb{AG}}(R)$. $\square$\\

 Proposition 1.8,  is not true in general, if we replace
$Spec(R)$ with $Max(R)$. For instance, if
$R=\Bbb{Z}\times\Bbb{Z}$, though  $R$ has $ACC$ on its ideals and
$\Bbb{A}(R)^*\neq\emptyset$, we have
$Max(R)\cap\Bbb{A}(R)=\emptyset$. In  the following result we
characterize all non-domain  rings $R$ for which
 $Max(R)\cap\Bbb{A}(R)\neq\emptyset$.\\

\noindent{\bf Proposition  1.12.} {\it Let $R$ be a non-domain
ring . Then $Max(R)\cap\Bbb{A}(R)^*\neq\emptyset$ if and only if
$soc(R)\neq (0)$.}\\

\noindent{\bf  Proof.} Suppose that $P\in Max(R)\cap\Bbb{A}(R)^*$.
Then  by Lemma 1.10, $P=Ann(x)$ for some $0\neq x\in R$. Since
$R/P\cong Rx$, $Rx$ is a minimal ideal of $R$, i.e., $soc(R)\neq
(0) $. Conversely, suppose that $soc(R)\neq (0)$ and $I$ is a
minimal ideal of $R$. Then $I=Rx$ for some $x\in I$. Since
$R/Ann(x)\cong Rx$ and  $Rx$ is a minimal ideal of $R$, $Ann(x)$ is a maximal ideal of $R$  and so $Max(R)\cap\Bbb{A}(R)^*\neq\emptyset$. $\square$\\

We conclude this section with the following  result,  that gives
us a characterization for  rings $R$ for which every nonzero
proper cyclic  ideal $I$  of $R$ is a  vertex of ${\Bbb{AG}}(R)$. The proof is trivial and left to reader. \\

\noindent{\bf Proposition  1.13.} {\it Let $R$ be a ring. Then
very nonzero proper cyclic ideal of $R$ is a vertex of
${\Bbb{AG}}(R)$ if and only if every element in $R$ is a unit or a
zero-divisor.}\\

\noindent{\bf 2. Connectivity of the annihilating-ideal graphs}\\

 By Anderson and Livingston [3, Theorem 2.3], for every
ring $R$, the zero divisor graph  $\Gamma(R)$ is a connected graph
and $diam(\Gamma(R))\leq  3$. Moreover, if $\Gamma(R)$ contains a
cycle, then gr( $\Gamma(R))\leq 4$ (see [16]). These facts later
were developed by Behboodi [8] for modules over a commutative
ring, by Redmond [17], for the undirected zero-divisor graph of a
non-commutative ring and  by Behboodi and Beyranvand [7] for the
strong zero divisor graphs of non-commutative rings.  Here we will
show that these facts are also true for the annihilating-ideal
graph of a   ring.

Let $S$ be a commutative multiplicative  semigroup with $0$
($0x=0$ for all $x\in S$). The zero-divisor graph of  $S$ (denoted
by $\Gamma(S)$) is a graph  whose vertices are the nonzero
zero-divisors of $S$,  with two distinct vertices $a$, $b$ joined
by an
 edge in case $ab=0$. In [11, Theorem 1.2 ], it is shown that $\Gamma(S)$ is always connected,
and $diam(\Gamma(S))\leq 3$.  Next, we use this to obtain the same result for the annihilating-ideal graph of a ring. \\

 \noindent{\bf Theorem 2.1.} {\it For every  ring $R$, the annihilating-ideal graph   ${\Bbb{AG}}(R)$ is  connected  and
$diam(\Bbb{AG}(R))\leq 3$. Moreover, if ${\Bbb{AG}}(R)$
contains a cycle, then $gr({\Bbb{AG}}(R))\leq 4$.}\\

 \noindent  {\bf  Proof.} Let $R$ be a    ring. Clearly,   the set
$S:=\Bbb{A}(R)$ is a commutative semigroup under multiplication
and also ${\Bbb{AG}}(R)=\Gamma(S)$. Thus by [11, Theorem 1.2 ],
${\Bbb{AG}}(R)$ is a connected graph and
$diam({\Bbb{AG}}(R))\leq 3$.\\
Now, suppose that ${\Bbb{AG}}(R)$ contains a cycle, and let
$C:=I_1-\hspace{-2mm}-\hspace{-2mm}-...
-\hspace{-2mm}-\hspace{-2mm}-I_n-\hspace{-2mm}-\hspace{-2mm}-I_1$
be a cycle with the least length. If $n\leq 4$, we are done.
Otherwise, we have  $I_1\cap I_4\neq (0)$. We need only consider 3 cases:\\
     \noindent Case 1:  $I_1\cap I_4=I_1$. Then
     $I_1I_3\subseteq I_4I_3=(0)$ and $I_{1}-\hspace{-2mm}-\hspace{-2mm}-I_{2}-\hspace{-2mm}-\hspace{-2mm}-I_{3}-\hspace{-2mm}-\hspace{-2mm}-I_{1}$
    is a cycle, a contradiction. The case $I_1\cap I_4=I_4$ is similar.\\
 \noindent Case 2:  $I_1\cap I_4=I_2$. Then
 $I_2\subseteq I_1$, $I_2I_n=(0)$ and so $I_2-\hspace{-2mm}-\hspace{-2mm}-...
-\hspace{-2mm}-\hspace{-2mm}-I_n-\hspace{-2mm}-\hspace{-2mm}-I_2$
is a cycle with length $n-1$, a contradiction. The case $I_1\cap I_4=I_3$ is similar.\\
    \noindent Case 3: $I_1\cap I_4\neq I_1,I_2,I_3,I_4$. Then  we have $I_{2}(I_{1}\cap I_{4})=(0)$ , $I_{3}(I_{1}\cap I_{4})=(0)$,
    and $I_{2}-\hspace{-2mm}-\hspace{-2mm}-(I_{1}\cap
    I_{4})-\hspace{-2mm}-\hspace{-2mm}-I_{3}-\hspace{-2mm}-\hspace{-2mm}-I_{2}$
    is a cycle, a contradiction. \\
    Thus $n\leq 4$,  i.e.,  $gr({\Bbb{AG}}(R))\leq 4$. $\square$\\

Next, we characterize all rings $R$ for which the graph
${\Bbb{AG}}(R)$ has  a vertex adjacent to every other vertex. Then
we apply this  to characterize  rings $R$ for which
the graph ${\Bbb{AG}}(R)$ is  complete or  star.\\

 \noindent{\bf Theorem 2.2.} Let $R$  be a
ring. Then there is a vertex of ${\Bbb{AG}}(R)$ which is adjacent
to every other vertex if and only if either $R=F\oplus D$, where
$F$ is a field and $D$ is an integral domain, or $Z(R)$ is an annihilator ideal.\\

 \noindent{\bf  Proof.} $(\Rightarrow)$. Suppose that
$Z(R)$  is not an annihilator ideal and the vertex $I_0\in
\Bbb{A}^*(R)$ is adjacent to every other vertex. Let $0\neq a\in
I_0$. Then $Ra$ is also adjacent to every other vertex. Now
$a\not\in Ann(a)=I$, for otherwise $Z(R)$  would be an annihilator
ideal.
 If $Rb$  is a nonzero  ideal of $R$ such that
$Rb\subsetneqq Ra$, then $Rb$ is also adjacent to every other
vertex and also $RbRb\subseteq RbRa=(0)$. Hence  for each $x\in
Z(R)$, $Rx$ is a vertex of ${\Bbb{AG}}(R)$ and   $(Rx)(Rb)=(0)$
i.e., $Z(R)=Ann(b)$, a contradiction. Thus $Ra$ is a minimal ideal
of $R$ with $(Ra)^2\neq (0) $. Thus by Brauer's Lemma (see [13,
10.22]), $Ra=Re$ for some idempotent $e\in R$. So $R=Re\oplus
R(1-e)$ and hence we may assume that $R=R_1\times R_2$ with
$R_1\times (0)$ adjacent to every other vertex. For each  $0\neq
c\in R_1$, $R_1c\times (0)$  is an annihilating-ideal of $R$. If
$R_1c\neq R_1$, then $(R_1\times (0))(R_1c\times
(0))=(0)\times(0)$, i.e., $cR_1=(0)$,  a contradiction. Thus $R_1$
must be  a field.  If $R_2$ is not an integral domain, then there
is a nonzero $b\in Z(R_2)$. Then $R_1\times R_2b$ is an
annihilating-ideal of $R$ which is not adjacent to $R_1\times
(0)$, a contradiction.
Thus $R_2$ must be an integral domain. \\
$(\Leftarrow)$.  If $R=F\oplus D$, where $F$ is a field and $D$ is
an integral domain, then $F\times (0)$ is adjacent to every other
vertex. If $Z(R)=Ann(x)$  for some nonzero $x\in R$, then $Rx$ is
adjacent to every other vertex. $\square$\\

Next, we characterize reduced rings $R$ for which the annihilating-ideal graph  ${\Bbb{AG}}(R)$ is a star graph.\\

 \noindent{\bf Corollary 2.3.} {\it Let $R$ be  a reduced
ring. Then  the following statements are equivalent.}\vspace{3mm}\\
(1) {\it There is a vertex of ${\Bbb{AG}}(R)$ which is adjacent
to every other vertex.}\vspace{2mm}\\
(2) {\it ${\Bbb{AG}}(R)$ is a star graph.}\vspace{2mm}\\
(3) {\it  $R=F\oplus D$, where
$F$ is a field and $D$ is an integral domain.}\\

\noindent {\bf  Proof.} $(1)\Rightarrow (3)$. Suppose there  is a
vertex of ${\Bbb{AG}}(R)$ which is adjacent to every other vertex.
If $Z(R)=Ann(x)$ for some $0\neq x\in R$, then $x^2=0$ and so
$x=0$, since $R$ is reduced. Thus by  Theorem 2.2, $R=F\oplus D$,
where
$F$ is a field and $D$ is an integral domain.\\
$(3)\Rightarrow (2)$. Let $R=F\oplus D$, where $F$ is a field and
$D$ is a domain. Then every nonzero ideal of $R$ is of the form
$F\oplus I$ or $(0)\oplus I$ where $I$ is a nonzero proper ideal
of $R$. By our hypothesis, and since ${\Bbb{AG}}(R)$ is connected,
we don't have any vertices of the form $F\oplus I $, such that
$I\neq (0)$. Also $F\oplus (0)$ is adjacent to every other vertex,
and since $D$ is a domain, non of the ideals of the form
$(0)\oplus I$ can be adjacent to each other.
So ${\Bbb{AG}}(R)$ is a star graph.\\
$(2)\Rightarrow (1)$ is clear.  $\square$\\

\noindent{\bf Corollary  2.4.} Let $R$  be a   Artinian ring. Then
there is a vertex of ${\Bbb{AG}}(R)$ which is adjacent to every
other vertex if and only if either $R=F_1\oplus F_2$, where
$F_1$, $F_2$ are fields, or $R$ is a local ring with nonzero maximal ideal.\\

\noindent{\bf  Proof.} $(\Rightarrow)$. Assume  there is a vertex
of ${\Bbb{AG}}(R)$ which is adjacent to every other vertex. Then
by Theorem 2.2,  either $R=F\oplus D$, where $F$ is a field and
$D$ is an integral domain, or $Z(R)$ is an annihilator ideal.
Since $R$ is Artinian, if $R=F\oplus D$, then $D$ is an Artinian
domain and so $D$ is also a field. Thus $R=F_1\oplus F_2$, where
$F_1$, $F_2$ are fields. If $Z(R)$ is  an annihilator ideal, then
$Z(R)$ is the set of all  non-unit elements of $R$ (since $R$ is
Artinian). It follows that every maximal
ideal of $R$ is contained in $Z(R)$, i.e.,  $R$ is  local with maximal ideal $M=Z(R)$.\\
$(\Leftarrow)$. If $R=F_1\oplus F_2$, where $F_1$, $F_2$ are
fields, then the graph ${\Bbb{AG}}(R)$  is a connected graph with
two vertices $F_1\times (0)$ and $(0)\times F_2$. Assume that $R$
is  an Artinian  local ring with nonzero maximal ideal $M$. Then
for each minimal ideal $I$ of $R$, $Ann(I)=M$, and hence $I$ is a
vertex of ${\Bbb{AG}}(R)$ which is adjacent to every other
vertex. $\square$\\

To have a better characterization  for an Artinian ring with star annihilating-ideal graph, we need the following Lemma.\\

 \noindent{\bf Lemma 2.5.} {\it Let $R$  be an Artinian  ring such that ${\Bbb{AG}}(R)$ is a star graph. Then
either $R=F_1\oplus F_2$, where $F_1$, $F_2$ are  fields, or $R$ is a local ring with nonzero maximal ideal  $M$ with  $M^4=(0)$.}\\

\noindent{\bf  Proof.} Let $R$  be an   Artinian ring such that
${\Bbb{AG}}(R)$ is a star graph. Then by Theorem 2.4, either
$R=F_1\oplus F_2$, where $F_1$, $F_2$ are  fields, or
 $R$ is a local ring with  nonzero maximal ideal  $M$. If $R=F_1\oplus F_2$, there is nothing to prove.
Suppose $R$ is  a local ring with   nonzero maximal ideal $M$.
Since $R$ is Artinian, there exists an integer $n\geq 1$ such that
$M^n=(0)$ and $M^{n-1}\neq (0)$. Clearly $M^{n-1}I\subseteq
M^{n-1}M=(0)$, for each ideal $I$ of $R$. Hence $M^{n-1}$ is
adjacent to every nonzero ideal $I$ of $R$. If $n>4$, then
 $M^2$ and $M^{n-2}$ will be adjacent, a contradiction.
 Thus $M(R)^4=(0)$. $\square$\\

\noindent{\bf Theorem 2.6.} {\it Let $R$  be an Artinian ring.
Then ${\Bbb{AG}}(R)$ is a star graph if and only if either
$R=F_1\oplus F_2$, where $F_1$, $F_2$ are fields,
or  $R$ is a local ring with nonzero maximal ideal $M$  and one of the following cases holds.}\vspace{3mm}\\
(i) {\it $M^2=(0)$ and $M$ is the only  nonzero proper ideal of $R$.}\vspace{2mm}\\
(ii) {\it $M^3=(0)$,  $M^2$ is the only minimal ideal of $R$ and
for every distinct \indent proper  ideals  $I_1$, $I_2$  of $R$
such
that $M^2\neq I_i$ $(i=1,2)$, $I_1I_2=M^2$.}\vspace{2mm}\\
(iii) {\it $M^4=(0)$, $M^3\neq (0)$ and
${\Bbb{A}^*}(R)=\{M,M^2,M^3\}$.}\\

  \noindent{\bf  Proof.} $(\Rightarrow)$. By Lemma 2.5,  either $R=F_1\oplus F_2$, where
$F_1$, $F_2$ are fields, or $R$ is a local ring with nonzero
maximal ideal $M$  such that $M^4=(0)$.We proceed by the following cases:\\
 \noindent Case 1:  $M^2=(0)$. It follows that
$M$ is the vertex adjacent to all other vertices.  Assume there
are more than  one nonzero proper ideals. Then there are nonzero
elements $x$, $y\in M$ with $xR\neq yR$. Since ${\Bbb{AG}}(R)$ is
a star graph and  $M^2=(0)$, either $xR=M$ or $yR=M$. Without loss
of generality we can assume that $yR=M$. Thus  $xR\subsetneqq
M=yR$ and hence $x=yr$ for some $r\in R$. Since $xR\neq yR$, $r$
is not a  unit element of $R$  and so $r\in M$. It follows that
$x=0$, a
contradiction.  Thus $M$ is the  only  nonzero proper ideal of $R$.\\
 \noindent Case 2: $M^3=(0)$ and $M^2\neq (0)$. Clearly $M^2$
is the vertex adjacent to all other vertices. Suppose that  $I$ is
a minimal ideal of $R$ such that $I\neq M^2$.  Since $Ann(I)=M$,
$I$ is also adjacent to all other vertices,  a contradiction. Thus
$M^2$ is the only minimal ideal of $R$. Now we assume that $I_1$,
$I_2$ are distinct proper  ideals  of $R$ such that $M^2\neq I_i$
($i=1,2$). Clearly, $I_1I_2\subseteq M^2$. On the other hand
$I_1I_2\neq (0)$ (since ${\Bbb{AG}}(R)$ is a star graph). Hence  $I_1I_2=M^2$. \\
\noindent Case 3: $M^4=(0)$ and $M^3\neq (0)$. Since
${\Bbb{AG}}(R)$ is a star graph and $R$ is local, then clearly the
"center" of the star must be a nonzero principal ideal. In the
case $M^3\neq M^4=(0)$, this ideal must be $M^3$. Since
$M^2\supsetneqq M^3$, there are elements $a$, $b\in M\setminus
M^2$ such that $ab\in M^2\setminus M^3$. Then $abM^2=(0)$, so the
star shaped assumption yields $M^2=abR$. If $a^2\in  M^3$, then
 $(Ra)(abR)=(0)$, contradicting the star shaped
assumption. Hence $M^2=a^2R$. Another application of the star
shaped assumption yields $a^3\neq 0$. So $M^3=a^3R$. For $c\in
M\setminus M^2$, $ca^2\in a^3R\setminus \{0\}$. From this point an
elementary argument puts $c\in aR$. Hence $M=aR$ and
${\Bbb{A}^*}(R)=\{M,M^2,M^3\}$.\\
 $(\Leftarrow)$ is clear. $\square$\\

 Now we are in position to characterize  rings with complete annihilating-ideal graphs.\\

\noindent{\bf Theorem  2.7.} {\it Let $R$  be a   ring. Then
${\Bbb{AG}}(R)$ is a complete graph if and only if $R$ is one of
the following three types of rings:}\vspace{3mm}\\
(1)  {\it  $R=F_1\oplus F_2$ where $F_1$, $F_2$ are fields,}\vspace{2mm}\\
(2) {\it  $Z(R)$ is an ideal of $R$ with  $Z(R)^2=(0)$, or}\vspace{2mm}\\
(3) {\it  $R$ is a local ring with exactly two nonzero proper
ideals $Z(R)$ and $Z(R)^2$.}\\

\noindent{\bf  Proof.} Assume that $R$ is a   ring for which
${\Bbb{AG}}(R)$ is a complete graph. Then, by  Theorem 2.2, either
$R=F\oplus D$, where $F$ is a field and $D$ is an integral domain,
or $Z(R)$ is an annihilator ideal. Suppose that $R=F\oplus D$,
where $F$  is a field and $D$  is an integral domain. If $D$ has a
nonzero proper ideal $I$, then $(0)\oplus D$ and $(0)\oplus I$ are
vertices of ${\Bbb{AG}}(R)$ which are not adjacent, a
contradiction. Thus $D$ doesn't have any nonzero proper ideal and
so it is a field. Let $Z(R)$ be an ideal of $R$. A trivial case is
$Z(R)^2=(0)(\neq Z(R))$, so assume not and let $x$, $y\in
Z(R)\setminus\{0\}$ be such that $xy\neq 0$. Then it must be that
$xR=yR$ with $x^2\neq  0 \neq  y^2$. If $z\in Z(R)\setminus\{0\}$
is such that $zR\neq xR$, then $xz=0$ and either $(x+z)R\neq zR$
or $(x+z)R\neq xR$. If the latter, $0=x(x+z)=x^2+xz=x^2$, a
contradiction. Thus $(x+z)R=xR$, $z\in xR$ and $z^2=0$ (using
$0=z(x+z)=z^2)$. Therefore $xR=Z(R)$ with $x^3=0(\neq x^2)$, and
$x^2R=Z(R)^2$ is the only other nonzero ideal of $R$.  The converse is clear. $\square$\\

\noindent{\bf Remark  2.8.} Let $R$ be a ring. In [3, Theorem
2.8], it is shown the zero divisor graph $\Gamma(R)$ is a complete
graph if and only if either $R=\Bbb{Z}_2\oplus \Bbb{Z}_2$ or
$Z(R)$ is an ideal of $R$ and $Z(R)^2=(0)$.
 Thus by  above theorem if $\Gamma(R)$ is a complete  graph, then  ${\Bbb{AG}}(R)$ is also a complete  graph, but the converse is not true even if $Z(R)$
 is an ideal of $R$. For example, for the ring $\Bbb{Z}_{p^3}$, where $p$ is a prime number, ${\Bbb{AG}}(R)$ is a complete  graph but
 $Z(R)^2\neq (0)$ and so $\Gamma(R)$ is not  a complete  graph.\\

To close above discussions, we ascertain all rings their
annihilating-ideal graphs have less than four vertices.\\

 \noindent{\bf Corollary  2.9.} {\it Let $R$ be a   ring.
 Then}\vspace{3mm}\\
 (a) {\it ${\Bbb{AG}}(R)$  is a graph with one vertex if
 and only if $R$ has only one nonzero \indent ideal.}\vspace{2mm}\\
 (b) {\it ${\Bbb{AG}}(R)$  is a graph with two  vertices if
 and only if $R=F_1\oplus F_2$, where \indent $F_1$,  $F_2$   are fields or  $R$ is an Artinian  local ring with
  exactly two nonzero \indent  proper   ideals $Z(R)$ and $Z(R)^2$.}\vspace{2mm}\\
 (c) {\it ${\Bbb{AG}}(R)$  is a graph with three   vertices if
 and only if $R$ is an Artinian   \indent local  ring with
  exactly three  nonzero  proper ideals $Z(R)$, $Z(R)^2$ and  \indent $Z(R)^3$.}\\

  \noindent{\bf  Proof.} (a) By Theorem 1.4 is clear.\\
  (b) Suppose that ${\Bbb{AG}}(R)$  is a graph with two  vertices. Since ${\Bbb{AG}}(R)$ is connected, then
  ${\Bbb{AG}}(R)$ is a complete (or star) graph. Thus by Theorem 2.6(i) and Theorem  2.7,
    $R$ should be of the form $F_1\oplus F_2$, or it should be a local
 ring with exactly two ideals $Z(R)$ and $Z(R)^2$. The converse is clear.\\
 (c) Suppose that ${\Bbb{AG}}(R)$  is a graph with three
   vertices. Then by Theorem 1.4, $R$ is an Artinian ring with  exactly three nonzero  proper ideals.
  Since ${\Bbb{AG}}(R)$ is connected,
  either it is a complete or a star graph. If ${\Bbb{AG}}(R)$ is a  complete graph, Theorem  2.7,
  implies that  $Z(R)$ is an ideal of $R$ with
$Z(R)^2=(0)$.  If  $Z(R)=pR$ for some $p\in R$, then for each
$z\in Z(R)$ such that $zR\neq pR$, $z=pr$ for some non-unit
element $R$. Since $Z(R)^2=(0)$, $z=0$ and so  $Z(R)$ is the only
nonzero  proper ideal of $R$, a contradiction. Thus $Z(R)$ is not
principal. Let $I$ and $J$ be  two nonzero  proper ideals  of $R$
such that $I, J\subsetneqq Z(R)$ and $I\neq J$. Then for each
$0\neq z\in Z(R)$, either $Rz=I$ or $Rz=J$. It follows that
$Z(R)=I\bigcup J$. Since $Z(R)$ is an ideal, $Z(R)=I$ or $Z(R)=J$,
a contradiction.  Thus ${\Bbb{AG}}(R)$ is a star graph. By Theorem
2.6,  either $Z(R)^3=(0)$, $Z(R)^2\neq (0)$  and $R$ has
 exactly one  nonzero proper ideal $I$ such that $Z(R)^2\subsetneqq I\subsetneqq
 Z(R)$ or  $R$ has exactly three  nonzero proper ideals
$Z(R)$, $Z(R)^2$ and $Z(R)^3$. In the first case, if $Z(R)$ is not
principal, then for each $0\neq z\in Z(R)$, either $zR=I$ or
$zR=Z(R)^2$, i.e,  $Z(R)=I\bigcup Z(R)^2$,
 a contradiction.  Thus $R$ is local with principal maximal ideal $Z(R)=pR$ such that $p^2\neq p^3=0$.
 Let $z\in Z(R)$ such that $zR\neq Z(R)$. Then $z=pr$ for some non-unit element $r\in R$ (i.e, $r\in Z(R)$).
 Therefore   $zp^2=0$, and since ${\Bbb{AG}}(R)$ is a star graph,   $zR=p^2R$. Thus
  $p^2R$ is the only other nonzero  proper ideal of $R$, a contradiction. Thus $R$ has exactly three  nonzero proper ideals
$Z(R)$, $Z(R)^2$ and $Z(R)^3$. The converse is clear. $\square$\\

A natural question one may ask is ``what will happen if all
vertices of ${\Bbb{AG}}(R)$ are prime (resp., maximal) ideals of
$R$". That is, what we answered below, in Theorem 2.10. In fact,
we characterize all rings satisfy this conditions. \\

\noindent{\bf Theorem  2.10.} {\it Let $R$  be a
 ring such that $R$ is not a domain. If  every vertex of ${\Bbb{AG}}(R)$  is a  prime ideal of $R$
 $($i.e. ${\Bbb{A}^*}(R)\subseteq Spec(R)$$)$, then  either either $R=F_1\oplus F_2$ for a pair of  fields $F_1$ and $F_2$ or
   $R$ has only one nonzero proper  ideal.}\\

 \noindent{\bf Proof.} Let $R$  be a
 ring such that $R$ is not a domain, and suppose that  ${\Bbb{A}^*}(R)\subseteq Spec(R)$$)$.
  First, if $z$ is a nonzero zero divisor, then $Q=zR$
is a prime ideal. If $Q=Q^2$, then can assume $z=z^2$. The ideal
$N=(1-z)R$ must be a prime as well. Neither $Q$ nor $N$ can have
proper subideals. Hence in this case $R=F_1\oplus F_2$ for fields
$F_1$ and $F_2$. The other possibility is that $x^2=0$ for all
nonzero $x\in Z(R)$. As $xR$ has to be prime, $xR=Z(R)$ must be
the only nonzero proper ideal. [For $r\in  R\setminus Z(R)$,
$rxR=xR$ which leads to $r$ being a unit.]. $\square$\\

The following is now immediate.\\

 \noindent{\bf Corollary  2.11.}
{\it Let $R$  be a
 ring such that it is not a domain. Then  the following statements are
equivalent.}\vspace{3mm}\\
(1) {\it ${\Bbb{A}}(R)^*\subseteq Max(R)$, i.e., every vertex of ${\Bbb{AG}}(R)$  is a  maximal ideal of $R$.}\vspace{2mm}\\
(2) {\it ${\Bbb{A}}(R)^*=Max(R)$.}\vspace{2mm}\\
(3) {\it ${\Bbb{A}}(R)^*=Spec(R)$.}\vspace{2mm}\\
(4) {\it ${\Bbb{A}}(R)^*\subseteq Spec(R)$.}\vspace{2mm}\\
(5) {\it Either  $R=F_1\oplus F_2$, where $F_1$, $F_2$ are fields,
or  $R$ has only one nonzero  \indent  proper ideal.}\\

  \noindent {\bf Note:} In Part II we shall continue the study of
  this contraction.\\\\

\noindent{\bf  Acknowledgments}\vspace{2mm}

   \noindent{ This work was partially supported by the Center of Excellence of Algebraic
Methods and Application of Isfahan University of Technology. }\\

\noindent{\bf References}\vspace{2mm}\\
{\footnotesize
 \noindent$[1]$ S. Akbari, A.
Mohammadian,  Zero-divisor graphs of
non-commutative   rings, J. Algebra \indent 296  (2006)  462-479.\\
$[2]$ D. D. Anderson, M. Naseer,  Beck's coloring of a commutative
  ring, J. Algebra 159 \indent (1993)  500-514.\\
$[3]$ D. F. Anderson, P. S. Livingston,  The zero-divisor graph of
a  commutative  ring, J. \indent Algebra  217 (1999) 434-447.\\
$[4]$ D. F. Anderson, A. Badawi, The total graph of a commutative
ring, J.  Algebra in \indent press.\\
 $[5]$ M. Axtell, J. Coykendall, J.
Stickles, Zero-divisor graphs of polynomials and power \indent
series over   rings, Comm. Algebra,  33 (2005) 2043-2050.\\
$[6]$ I. Beck,  Coloring of commutative   rings, J. Algebra 116
(1988) 208-226. \\
$[7]$ M. Behboodi, R. Beyranvand, Strong zero-divisors graphs  of
non-commutative rings, \indent Int. J.  Algebra 2(1) (2008), 25-44.\\
  $[5]$ M. Behboodi, Zero divisor graphs of modules over a commutative   rings, J. Commuta- \indent tive Algebra, to appear. \\
$[6]$ A. Cannon, K. Neuerburg, S.P. Redmond, Zero-divisor graphs
of nearrings and  semi- \indent groups,   in: H. Kiechle, A.
Kreuzer,
M.J. Thomsen (Eds.), Nearrings and Nearfields, \indent Springer,  Dordrecht, The Netherlands, 2005, pp. 189-200.\\
$[10]$ F. DeMeyer, L. DeMeyer,  Zero divisor graphs of semigroups,
J. Algebra 283 (2005) \indent 190-198.\\
$[11]$ F. DeMeyer, T. McKenzie and K. Schneider,  The zero-divisor
graph of a commutative \indent semigroup,  Semigroup Forum 65 (2)
(2002), 206-214.\\
 $[12]$ I. Kaplansky, ''Commutative Rings,'' rev. ed., Univ. of Chicago Press, Chicago, 1974.\\
 $[13]$ T. Y. Lam,  A first course in non-commutative   rings,
Springer-Verlag New  York, Inc \indent 1991.\\
$[14]$ T. G. Lucas,  The diameter of a zero divisor graph. J.
Algebra 301 (2006) 174-193.\\
 $[15]$ H. R. Maimani, M. Salimi, A. Sattari,  S. Yassemi, Comaximal graph of  commutative  \indent rings,
    J. Algebra 319 (2008) 1801-1808.\\
 $[16]$ S. B. Mulay,  Cycles and symmetries of zero-divisors, Comm.
Algebra, 30 (2002) 3533- \indent 3558.\\
  $[17]$ S.P. Redmond,  The
zero-divisor graph of a non-commutative   ring, Internat.  J.
Com-   \indent  mutative Rings, 1 (4) (2002) 203-211.\\
$[18]$ P. K. Sharma, S. M. Bhatwadekar, A note on graphical
representation of rings, J. \indent Algebra  176 (1995) 124-127.\\
 $[19]$  R. Y. Sharp, Steps in   algebra, Second edition, London Mathematical Society Student \indent Texts,  51.
  Cambridge University Press, Cambridge, 1990.
  \end{document}